\documentclass{amsart}

\usepackage[T1]{fontenc}
\usepackage{palatino}
\usepackage{fullpage}
\usepackage{amsmath,amssymb,amsthm,mathtools}
\usepackage{xcolor}
\usepackage{thm-restate}
\usepackage{microtype}
\usepackage[colorlinks=true,linkcolor=blue,citecolor=red,urlcolor=blue]{hyperref}

\newtheorem{theorem}{Theorem}[section]
\newtheorem{corollary}[theorem]{Corollary}
\newtheorem{lemma}[theorem]{Lemma}
\newtheorem{proposition}[theorem]{Proposition}

\theoremstyle{definition}
\newtheorem{definition}[theorem]{Definition}
\newtheorem{example}[theorem]{Example}
\newtheorem{remark}[theorem]{Remark}
\newtheorem{question}{Question}

\numberwithin{equation}{section}

\newcommand{\N}{\mathbb N}
\newcommand{\Z}{\mathbb Z}
\newcommand{\R}{\mathbb R}
\newcommand{\X}{\mathcal X}
\newcommand{\Rcal}{\mathcal R}
\newcommand{\Bcal}{\mathcal B}
\DeclarePairedDelimiter{\abs}{\lvert}{\rvert}

\title{Gaps in the Most Frequent Sizes of Integer Sumsets}

\author{Saad A. Awan}
\address[Saad A. Awan]{Department of Mathematics, University of Kansas, Lawrence, KS 66045}
\email{\href{mailto:saad_awan@ku.edu}{saad\_awan@ku.edu}}

\author{Steven J. Miller}
\address[Steven J. Miller]{Department of Mathematics, Williams College, Williamstown, MA 01267}
\email{\href{mailto:sjm1@williams.edu}{sjm1@williams.edu}, \href{mailto:Steven.Miller.MC.96@aya.yale.edu}{Steven.Miller.MC.96@aya.yale.edu}}

\author{Linh Nguyen}
\address[Linh Nguyen]{Department of Mathematics, Gettysburg College, Gettysburg, PA 17325}
\email{\href{mailto:nguyli01@gettysburg.edu}{nguyli01@gettysburg.edu}}

\author{Steven Senger}
\address[Steven Senger]{Department of Mathematics, Missouri State University, Springfield, MO 65809}
\email{\href{mailto:stevensenger@missouristate.edu}{stevensenger@missouristate.edu}}

\author{Elvis Wong}
\address[Elvis Wong]{Department of Mathematics, The University of Manchester, Manchester, M13 9PL}
\email{\href{mailto:yik.wong@student.manchester.ac.uk}{yik.wong@student.manchester.ac.uk}}

\begin{document}

\begin{abstract}
Nathanson observed that, for four-element sets $A$ of integers, the most frequent sizes of the iterated sumset $hA$ are separated by triangular numbers. Senger proved this and suggested an analogue for larger sets. We establish this pattern for $k$-element subsets of $\{1,\ldots,q\}$, for every fixed $k\geq3$, where $h\geq2$ is either fixed or grows slower than $q^{1/(2k-1)}$. As $q\to\infty$, almost every such set has the largest possible $h$-fold sumset, and almost every exception has one of $h-1$ explicit smaller sizes. The differences between consecutive sizes among these $h$ values are consecutive $(k-2)$-dimensional simplex numbers.

\end{abstract}

\maketitle

\section{Introduction}
Let $A\subset \mathbb{Z}$ be a finite subset.
The \textbf{sumset} of $A$ is
\[
A+A \ = \ \{a+a':a,a'\in A\} \ = \ 2A
\]
and its \textbf{$h$-fold sumset} is
\[
hA \ = \ (h-1)A+A \ = \ \{a+a':a\in A,\ a'\in(h-1)A\},
\]
where $1A=A$ and repetition is allowed.

The size of $hA$ has been studied in two complementary regimes. In the first, the set $A$ is fixed and $h\to\infty$. Nathanson~\cite{Nathanson1972} proved that if $A$ is a finite set of integers with $\min A=0$ and $\gcd(A)=1$, then for all sufficiently large $h$, the integers in $[0..h\max A]$ that are missing from $hA$ form two finite sets, one near each endpoint of this interval, which do not depend on $h$ up to translation. In particular, $\abs{hA}$ is eventually a linear function of $h$. Khovanskii~\cite{Khovanskii1992} proved that for a finite subset $A$ of $\Z^d$, $\abs{hA}$ agrees with a polynomial in $h$ for all sufficiently large $h$. Nathanson and Ruzsa~\cite{NathansonRuzsa2002} gave an elementary lattice-point proof of eventual polynomial growth for finite subsets of abelian semigroups. For finite subsets of $\Z^d$, Granville, Shakan and Walker~\cite{GranvilleShakanWalker2023} gave effective bounds on how large $h$ must be. A finite set $B$ of integers is \textbf{sum-dominant} if $\abs{B+B}>\abs{B-B}$. In the study of such sets, Iyer, Lazarev, Miller and Zhang~\cite{ILMZ2012} determined the limiting behavior of $hA$ for an arbitrary finite set $A$ as $h\to\infty$, together with a bound on how large $h$ must be for this behavior to appear. They also proved that no set $A$ has the property that $A,2A,\ldots,hA$ are all sum-dominant for every $h$.

In the second regime, which is the subject of this paper, $h$ is fixed or grows slowly with $q$, and $A$ is a $k$-element subset of $\{1,\dots, q\}$ with $q\to\infty$. A $k$-element set $A$ is a $B_h$-set when $\abs{hA}$ attains the universal upper bound $M_{h,k}$ (see Definition~\ref{def:Bh-collision}). For $h=2$ these are the Sidon sets, studied by Erd\H{o}s and Tur\'an~\cite{ErdosTuran1941}, and Bose and Chowla~\cite{BoseChowla1962} constructed large $B_h$-sets for every $h$. We refer to O'Bryant's survey~\cite{OBryant2004} for the extensive literature on these sets. Nathanson~\cite{Nathanson2004} proved that almost all $k$-element subsets of $[q]$ are Sidon sets, even when $k$ grows with $q$ at a suitable rate. Random subsets of intervals are also central to the study of sum-dominant sets. Martin and O'Bryant~\cite{MartinOBryant2007} proved that a positive proportion of subsets of $\{0,\ldots,n\}$ are sum-dominant, and Hegarty and Miller~\cite{HegartyMiller2009} proved that when the elements of $[q]$ are chosen independently with probability $p(q)$ satisfying $p(q)\to0$ and $qp(q)\to\infty$, almost all sets have more differences than sums.

One particular approach has been to study the cardinalities of $\lvert hA \rvert$ for various values of $h$. Nathanson~\cite[Problem 5]{NathansonTriangular} first studied the most frequent sizes of $\lvert hA \rvert$ for $\lvert A \rvert = 4$ and observed triangular gaps between them.
O'Bryant proved a formula for the initial sequence of sumset sizes
in terms of the successive $L^1$-minima of a lattice and gave a heuristic
explanation of the observed frequencies \cite{OBryantTriangular}. The fourth listed author studied this problem further, and proved
this phenomenon to be true for $k=4$
\cite[Thm.~2]{Senger}. He also suggested that for $k>4$, analogous
``champagnerpyramide'' numbers would describe the sumset-size deficits
\cite[Sec.~4]{Senger}. These are $(k-1)$-dimensional simplex numbers,
whose consecutive differences are $(k-2)$-dimensional simplex numbers.
We prove that this gap pattern holds for all fixed $k\geq3$ by generalizing the fourth author's constructions, as we now explain.

For $h\geq0$, let
\[
\X_{h,k}
\ = \
\left\{\mathbf x=(x_1,\ldots,x_k)\in\Z_{\geq0}^k:
\sum_{i=1}^k x_i=h\right\}.
\]
We call elements of $\mathcal{X}_{h,k}$ \textbf{sum vectors} since by representing elements of $A = \{a_1 < a_2 < \cdots < a_k\}$ as $\vec{A} = (a_1, a_2, \dots, a_k) \in \mathbb{Z}_{\ge 0}^{k}$ we can produce any element of $hA$. We may redefine $hA$ as
\[
hA  \ = \  \{\mathbf{x} \cdot \vec{A}: \mathbf{x} \in \mathcal{X}_{h,k}\}.
\]

Let
\[
M_{h,k} \ = \ \abs{\X_{h,k}} \ = \ \binom{h+k-1}{k-1}, \quad M_{-1,k} \ = \ 0
\]
be the maximum possible size of $hA$. We are ready to state our main result.

\begin{restatable}{theorem}{maintheorem}\label{thm:main}
Let $k\geq3$ be fixed, and let $q\to\infty$ and $h=h(q)\geq2$ with $h=o(q^{1/k})$. For
$0\leq\ell<h$, let
\[
P_\ell \ = \ M_{h,k}-M_{\ell-1,k}.
\]
Then the following statements hold.
\begin{enumerate}
\item
\[
\#\left\{
A\in\binom{[q]}k:
\abs{hA}=M_{h,k}
\right\}
\ = \
\Theta(q^k).
\]

\item
\[
\#\left\{A\in\binom{[q]}k:
\abs{hA}\notin\{P_0,\ldots,P_{h-1}\}\right\}
\ = \ O\!\left(h^{2k-2}q^{k-2}\right).
\]

\item
\[
\#\left\{
A\in\binom{[q]}k:
\abs{hA}=P_\ell
\right\}
\ = \
\Omega\!\left(
\frac{q^{k-1}}{h-\ell+1}
\right).
\]
\end{enumerate}
\end{restatable}

Accordingly, we show that the natural gaps between the most frequent sizes of $\lvert hA\rvert$ are of the form
\[
\binom{\ell+k-2}{k-2},
\]
which specialize to the triangular numbers when $k=4$.

\subsection*{Formalization}
The results of this paper have been formalized in the Lean~4 proof assistant
with the Mathlib library~\cite{AwanLean}. Part of the Lean code was written by
Awan, and part was generated by Aristotle, the automated theorem prover
developed by Harmonic~\cite{Aristotle}. Every proof, including the generated
ones, is checked by the Lean kernel, with the one exception noted at the end of
this paragraph. The formalization contains no \texttt{sorry}. It proves
Lemmas~\ref{lem:equal-sum-criterion}, \ref{lem:rank-lemma},
\ref{lem:gap-size-lemma}, \ref{lem:disjoint-support} and~\ref{lem:Bstar-count},
Propositions~\ref{prop:coefficient-vectors-lemma}, \ref{prop:T-dichotomy}
and~\ref{prop:S-interval}, Theorem~\ref{thm:main} for fixed $h$ together with explicit versions of
its three bounds, which hold once $q$ exceeds an explicit polynomial in $h$
and imply Theorem~\ref{thm:main} as stated,
Corollaries~\ref{cor:gap-pattern} and~\ref{cor:proportion},
Remark~\ref{rem:k3-limit}, and Examples~\ref{ex:E}, \ref{ex:E-relation}
and~\ref{ex:E-size}. In particular, it contains a proof of
Lemma~\ref{lem:disjoint-support}, for which we only cite Senger's argument, and
it proves Lemma~\ref{lem:gap-size-lemma} directly for every $k$, without using
O'Bryant's theorem. For Remark~\ref{rem:explicit}, the formalization proves the
conclusion of Corollary~\ref{cor:gap-pattern} for every
$q>70\cdot2^k\,k!\,h\,(2h+1)^{2k-2}$, which is $70$ times the threshold stated
there, and it checks the computation for $(h,k)=(3,4)$ and $16\leq q\leq50$ with
\texttt{native\_decide}, which relies on Lean's compiler in addition to its
kernel.

\subsection*{Acknowledgments}
We thank the organizers of the Polymath Jr. research group whose efforts brought about this collaboration. \textbf{This work was supported by the National Science Foundation under Grant DMS2341670.}
\section{Background}
Throughout we assume familiarity with standard definitions and notation of additive and algebraic combinatorics; see, e.g.,~\cite{Tao_Vu_2006, Stanley_2011}.

For integers $u\leq v$ and $q\geq1$, we write
\[
[u..v] \ = \ \{n\in\Z:u\leq n\leq v\},
\qquad [q] \ = \ [1..q].
\]

\begin{definition}\label{def:subsets}
Let $q\geq k\geq1$. The set of $k$-element subsets of $[q]$ is
\[
\binom{[q]}k \ = \ \{A\subseteq[q]:\abs{A} \ = \ k\}.
\]
For $A\in\binom{[q]}k$, we always list the elements of $A$ in increasing order, $a_1<a_2<\cdots<a_k$, and we call
\[
\vec{A} \ = \ (a_1,a_2,\ldots,a_k)\in[q]^k
\]
the \textbf{set vector} of $A$.
\end{definition}

The map $A\mapsto\vec{A}$ is injective, so we may count subsets of $[q]$ by counting their set vectors in $[q]^k$.

\begin{definition}\label{def:Bh-collision}
Let $A\in\binom{[q]}k$ and $h\geq1$. The set $A$ is a \textbf{$B_h$-set} if $\abs{hA}=M_{h,k}$. Following~\cite{Senger}, an element $n\in hA$ is a \textbf{collision} in $hA$ if $n$ can be written in at least two different ways as a sum of $h$ elements of $A$ listed in nondecreasing order.
\end{definition}

For instance, $a_1+a_2$ and $a_2+a_1$ are not two different ways of writing an element of $2A$, since only $a_1+a_2$ is listed in nondecreasing order. A sum of $h$ elements of $A$ listed in nondecreasing order is determined by the number of times each $a_i$ occurs in it, that is, by a sum vector in $\X_{h,k}$. Hence $n$ is a collision in $hA$ exactly when $n=\mathbf x\cdot\vec{A}=\mathbf y\cdot\vec{A}$ for distinct $\mathbf x,\mathbf y\in\X_{h,k}$. Consequently, $A$ is a $B_h$-set if and only if $hA$ has no collisions, that is, if and only if the map $\mathbf x\mapsto\mathbf x\cdot\vec{A}$ is injective on $\X_{h,k}$.

\begin{example}\label{ex:E}
Let $q\geq3$ and $E=\{1,2,3\}\in\binom{[q]}3$, so that $\vec{E}=(1,2,3)$. The integer $8$ is a collision in $4E$, because
\[
1+2+2+3 \ = \ 1+1+3+3 \ = \ 2+2+2+2 \ = \ 8.
\]
In terms of sum vectors in $\X_{4,3}$, this reads
\[
(1,2,1)\cdot\vec{E} \ = \ (2,0,2)\cdot\vec{E} \ = \ (0,4,0)\cdot\vec{E} \ = \ 8.
\]
Moreover, $4E=[4..12]$, so $\abs{4E}=9<15=M_{4,3}$ and $E$ is not a $B_4$-set.
\end{example}

\begin{definition}\label{def:asymptotic}
Let $f:\N\to\R$ and $g:\N\to\R_{>0}$. We write $f(q)=O(g(q))$ if there exist $C>0$ and $q_0\in\N$ such that $\abs{f(q)}\leq Cg(q)$ for all $q\geq q_0$, and $f(q)=\Omega(g(q))$ if there exist $c>0$ and $q_0\in\N$ such that $f(q)\geq cg(q)$ for all $q\geq q_0$. We write $f(q)=\Theta(g(q))$ if $f(q)=O(g(q))$ and $f(q)=\Omega(g(q))$.
\end{definition}

From here on, all asymptotic notation is with respect to $q\to\infty$, and $k\geq3$ is fixed. We treat $h=h(q)\geq2$ as a second asymptotic parameter, which may grow with $q$. The constants $C$ and $c$ depend only on $k$, and $q_0$ depends only on $k$ and on the function $h$. Each statement records the growth condition on $h$ that it needs. In Section~\ref{sec:additional}, $h$ is fixed, and the constants there may also depend on $h$.

\section{Preliminary Lemmas}

Throughout this section, $A\in\binom{[q]}k$ and $\vec{A}$ is its set vector. Collisions in $hA$ correspond to integer vectors $\mathbf d$ with coordinate sum zero and $\mathbf d\cdot\vec{A}=0$. We normalize these vectors as follows.

\begin{definition}\label{def:coefficient-vectors}
Let $\mathbf d=(d_1,\ldots,d_k)\in\Z^k$ be a nonzero vector with $\sum_i d_i=0$. Let $\mathbf d^+=(d_1^+,\ldots,d_k^+)$ and $\mathbf d^-=(d_1^-,\ldots,d_k^-)$, where
\[
d_i^+ \ = \ \max\{d_i,0\},
\qquad
d_i^- \ = \ \max\{-d_i,0\},
\]
so that $\mathbf d=\mathbf d^+-\mathbf d^-$. The \textbf{weight} of $\mathbf d$ is
\[
w(\mathbf d) \ = \ \sum_i d_i^+ \ = \ \sum_i d_i^-.
\]
The vector $\mathbf d$ is \textbf{primitive} if $\gcd(\abs{d_1},\ldots,\abs{d_k})=1$. Let $j(\mathbf d)$ be the least index $i$ with $d_i\neq0$.

The set of \textbf{coefficient vectors} of weight at most $h$ is
\begin{equation}\label{eq:collision-vector-set}
  \Rcal_{h,k}
  = \left\{ \mathbf d\in\Z^k\setminus\{\mathbf 0\} :
    \begin{aligned}
      &d_1+\dots+d_k = 0,                  && w(\mathbf d)\le h,\\
      &\gcd(\abs{d_1},\dots,\abs{d_k}) = 1, && d_{j(\mathbf d)} > 0
    \end{aligned}
  \right\}.
\end{equation}
The condition $d_{j(\mathbf d)}>0$ selects one vector from each pair $\{\mathbf d,-\mathbf d\}$.
\end{definition}

\begin{lemma}\label{lem:equal-sum-criterion}
There exist distinct $\mathbf x,\mathbf y\in\X_{h,k}$ such that
$\mathbf x\cdot\vec{A}=\mathbf y\cdot\vec{A}$ if and only if there
exists $\mathbf d\in\Rcal_{h,k}$ such that
$\mathbf d\cdot\vec{A}=0$.
\end{lemma}

\begin{proof}
Suppose that distinct $\mathbf x,\mathbf y\in\X_{h,k}$ satisfy
$\mathbf x\cdot\vec{A}=\mathbf y\cdot\vec{A}$. Let $\mathbf z=\mathbf x-\mathbf y$
and $g=\gcd(\abs{z_1},\ldots,\abs{z_k})$. Then
\begin{align*}
\sum_{i=1}^k z_i
& \ = \ \sum_{i=1}^k x_i-\sum_{i=1}^k y_i \ = \ h-h \ = \ 0,\\
\mathbf z\cdot\vec{A}
& \ = \ \mathbf x\cdot\vec{A}-\mathbf y\cdot\vec{A} \ = \ 0,\\
w(\mathbf z)
& \ = \ \sum_{i=1}^k\max\{x_i-y_i,0\}
  \ \leq\ \sum_{i=1}^k x_i \ = \ h.
\end{align*}
Choose $\varepsilon\in\{1,-1\}$ so that
$\varepsilon z_{j(\mathbf z)}>0$, and let $\mathbf d=\varepsilon\mathbf z/g$.
Then
\begin{align*}
\sum_{i=1}^k d_i
& \ = \ \frac{\varepsilon}{g}\sum_{i=1}^k z_i \ = \ 0,\\
\gcd(\abs{d_1},\ldots,\abs{d_k})
& \ = \ \frac{g}{g} \ = \ 1,\\
w(\mathbf d)
& \ = \ \frac{w(\mathbf z)}{g} \ \leq\ h,\\
d_{j(\mathbf d)}
& \ = \ \frac{\abs{z_{j(\mathbf z)}}}{g} \ > \ 0,\\
\mathbf d\cdot\vec{A}
& \ = \ \frac{\varepsilon}{g}\mathbf z\cdot\vec{A} \ = \ 0.
\end{align*}
Thus $\mathbf d\in\Rcal_{h,k}$.

Conversely, suppose that $\mathbf d\in\Rcal_{h,k}$ and
$\mathbf d\cdot\vec{A}=0$. Let $s=w(\mathbf d)$, choose
$\mathbf t\in\X_{h-s,k}$, and set
$\mathbf x=\mathbf d^++\mathbf t$ and $\mathbf y=\mathbf d^-+\mathbf t$.
Then
\begin{align*}
\sum_{i=1}^k x_i
& \ = \ \sum_{i=1}^k d_i^++\sum_{i=1}^k t_i
  \ = \ s+(h-s) \ = \ h,\\
\sum_{i=1}^k y_i
& \ = \ \sum_{i=1}^k d_i^-+\sum_{i=1}^k t_i
  \ = \ s+(h-s) \ = \ h,\\
\mathbf x-\mathbf y
& \ = \ \mathbf d^+-\mathbf d^- \ = \ \mathbf d \ \neq\ \mathbf 0,\\
\mathbf x\cdot\vec{A}-\mathbf y\cdot\vec{A}
& \ = \ (\mathbf x-\mathbf y)\cdot\vec{A}
  \ = \ \mathbf d\cdot\vec{A} \ = \ 0.
\end{align*}
Since $\mathbf x$ and $\mathbf y$ are nonnegative integer vectors, these
equalities give distinct $\mathbf x,\mathbf y\in\X_{h,k}$ with
$\mathbf x\cdot\vec{A}=\mathbf y\cdot\vec{A}$.
\end{proof}

\begin{example}\label{ex:E-relation}
In Example~\ref{ex:E}, the equality $(1,2,1)\cdot\vec{E}=(0,4,0)\cdot\vec{E}$ gives $(1,-2,1)\cdot\vec{E}=0$, and $(1,-2,1)\in\Rcal_{4,3}$. This vector encodes the relation $1+3=2+2$. In fact, for every $h\geq2$ it is the only vector in $\Rcal_{h,3}$ orthogonal to $\vec{E}$. Indeed, if $\mathbf d\cdot\vec{E}=0$ and $d_1+d_2+d_3=0$, then subtracting the second equation from the first gives $d_2+2d_3=0$, so $\mathbf d=d_3(1,-2,1)$.
\end{example}

\begin{lemma}\label{lem:rank-lemma}
Let $\mathbf d^{(1)},\ldots,\mathbf d^{(r)}\in\Z^k$ be linearly
independent. At most $q^{k-r}$ vectors $\mathbf v\in[q]^k$ obey
\[
\mathbf d^{(1)}\cdot\mathbf v
=\cdots=
\mathbf d^{(r)}\cdot\mathbf v=0.
\]
\end{lemma}

\begin{proof}
Let $D$ be the $r\times k$ matrix
\[
D=
\begin{pmatrix}
d_1^{(1)} & d_2^{(1)} & \cdots & d_k^{(1)}\\
d_1^{(2)} & d_2^{(2)} & \cdots & d_k^{(2)}\\
\vdots & \vdots & \ddots & \vdots\\
d_1^{(r)} & d_2^{(r)} & \cdots & d_k^{(r)}
\end{pmatrix}.
\]
Since its rows are linearly independent,
\[
\operatorname{rank}(D)=r.
\]
By equality of row rank and column rank, there exist
$r$ linearly independent columns indexed by
$1\le j_1<\cdots<j_r\le k$. Thus the $r\times r$
submatrix
\[
D'=
\begin{pmatrix}
d_{j_1}^{(1)} & \cdots & d_{j_r}^{(1)}\\
\vdots & \ddots & \vdots\\
d_{j_1}^{(r)} & \cdots & d_{j_r}^{(r)}
\end{pmatrix}
\]
satisfies
\[
\det(D')\neq 0
\implies
(D')^{-1}\text{ exists}.
\]

Let $J=\{j_1,\ldots,j_r\}$ and fix the $k-r$
coordinates $v_j$ with $j\notin J$. Writing
$\mathbf v'=(v_{j_1},\ldots,v_{j_r})^T$, we obtain
\[
D\mathbf v=\mathbf 0
\iff
D'\mathbf v'
=
-\begin{pmatrix}
\displaystyle\sum_{j\notin J}d_j^{(1)}v_j\\
\vdots\\
\displaystyle\sum_{j\notin J}d_j^{(r)}v_j
\end{pmatrix}
=\mathbf b
\]
and hence
\[
D'\mathbf v'=\mathbf b
\implies
\mathbf v'=(D')^{-1}\mathbf b.
\]
Thus each choice of the $k-r$ fixed coordinates
determines at most one vector $\mathbf v'$. Since
each coordinate has $q$ possible values, there
are at most
\[
q^{k-r}
\]
vectors $\mathbf v\in[q]^k$ satisfying
$D\mathbf v=\mathbf 0$.
\end{proof}

\begin{proposition}\label{prop:coefficient-vectors-lemma}
For $h\geq2$,
\[
\abs{\Rcal_{h,k}} \ \leq\ (2h+1)^{k-1} \ \leq\ 3^{k-1}h^{k-1}.
\]
Moreover, any two distinct vectors in $\Rcal_{h,k}$ are linearly
independent.
\end{proposition}

\begin{proof}
For $\mathbf d\in\Rcal_{h,k}$, we have
\[
-h\leq-w(\mathbf d)\leq d_i\leq w(\mathbf d)\leq h,
\qquad
d_k \ = \ -\sum_{i=1}^{k-1}d_i.
\]
The first $k-1$ coordinates therefore determine $\mathbf d$, so
\begin{align*}
\abs{\Rcal_{h,k}}
& \ = \ \#\{(d_1,\ldots,d_{k-1}):\mathbf d\in\Rcal_{h,k}\}\\
& \ \leq\ \abs{[-h..h]^{k-1}}
  \ = \ (2h+1)^{k-1}
  \ \leq\ 3^{k-1}h^{k-1}.
\end{align*}

Suppose that $\mathbf d,\mathbf e\in\Rcal_{h,k}$ are linearly dependent.
Since both are nonzero integer vectors, write
$\mathbf e=(a/b)\mathbf d$ with $a,b\in\Z$, $b>0$ and $\gcd(\abs{a},b)=1$.
Then
\begin{align*}
be_i=ad_i\quad\text{for every }i
& \ \implies\ b\mid d_i\quad\text{for every }i\\
& \ \implies\ b\mid\gcd(\abs{d_1},\ldots,\abs{d_k})=1\\
& \ \implies\ b=1.
\end{align*}
Consequently,
\[
1 \ = \ \gcd(\abs{e_1},\ldots,\abs{e_k})
  \ = \ \abs{a}\gcd(\abs{d_1},\ldots,\abs{d_k})
  \ = \ \abs{a}.
\]
Thus $a\in\{1,-1\}$. The sign convention in
\eqref{eq:collision-vector-set} gives
\[
0<e_{j(\mathbf e)}=a d_{j(\mathbf d)},
\qquad d_{j(\mathbf d)}>0
\quad\implies\quad a=1
\quad\implies\quad\mathbf e=\mathbf d.
\]
Hence distinct vectors in $\Rcal_{h,k}$ are linearly independent.
\end{proof}

\begin{lemma}[{cf.~\cite[Thm.~1]{OBryantTriangular}}]\label{lem:gap-size-lemma}
Suppose there is exactly one $\mathbf d\in\Rcal_{h,k}$ for which
$\mathbf d\cdot\vec{A}=0$. If $s=w(\mathbf d)$, then
\begin{equation}\label{eq:unique-vector-size}
\abs{hA} \ = \ M_{h,k}-M_{h-s,k}.
\end{equation}
\end{lemma}

\begin{proof}
For $k\geq4$, the lemma follows from O'Bryant's Theorem~1~\cite{OBryantTriangular},
because the hypothesis says that the lattice of all $\mathbf v\in\Z^k$ with
$v_1+\cdots+v_k=0$ and $\mathbf v\cdot\vec{A}=0$ has first $L^1$-minimum $2s$ and
second $L^1$-minimum greater than $2h$.

Now let $k=3$, and write
\[
\frac{a_2-a_1}{a_3-a_1} \ = \ \frac bc,
\qquad \gcd(b,c)=1,
\qquad 0<b<c.
\]
Let $\mathbf v\in\Z^3$ satisfy $v_1+v_2+v_3=0$ and $\mathbf v\cdot\vec{A}=0$. Then
\begin{align*}
0 \ = \ \mathbf v\cdot\vec{A}
& \ = \ (-v_2-v_3)a_1+v_2a_2+v_3a_3
  \ = \ v_2(a_2-a_1)+v_3(a_3-a_1)\\
& \ \implies\ v_2b(a_3-a_1) \ = \ v_2c(a_2-a_1) \ = \ -v_3c(a_3-a_1)
  \qquad(\text{because }c(a_2-a_1)=b(a_3-a_1))\\
& \ \implies\ v_2b \ = \ -v_3c
  \qquad(\text{because }a_3>a_1)\\
& \ \implies\ c\mid v_2
  \qquad(\text{because }\gcd(b,c)=1)\\
& \ \implies\ v_2=-mc,\quad v_3=mb,\quad v_1=-v_2-v_3=m(c-b)
  \qquad(m\in\Z)\\
& \ \implies\ \mathbf v \ = \ m\mathbf p,
  \qquad \mathbf p \ = \ (c-b,-c,b).
\end{align*}
Since $\mathbf d$ satisfies the same two conditions, $\mathbf d=m\mathbf p$ for some
$m\in\Z$, and
\begin{align*}
1 \ = \ \gcd(\abs{d_1},\abs{d_2},\abs{d_3})
& \ = \ \abs m\gcd(c-b,c,b) \ = \ \abs m\gcd(b,c) \ = \ \abs m\\
& \ \implies\ d_1 \ = \ m(c-b) \ \neq\ 0
  \ \implies\ j(\mathbf d) \ = \ 1\\
& \ \implies\ m(c-b) \ = \ d_{j(\mathbf d)} \ > \ 0
  \ \implies\ m \ = \ 1.
\end{align*}
Hence $\mathbf d=\mathbf p$, and
\begin{equation}\label{eq:k3-kernel}
\mathbf v\in\Z^3,\quad v_1+v_2+v_3=0,\quad \mathbf v\cdot\vec{A}=0
\quad\implies\quad \mathbf v\in\Z\mathbf d.
\end{equation}

For $\mathbf x,\mathbf y\in\X_{h,3}$, by \eqref{eq:k3-kernel} we get
\[
\mathbf x\cdot\vec{A}=\mathbf y\cdot\vec{A}
\quad\iff\quad
(\mathbf x-\mathbf y)\cdot\vec{A}=0
\quad\iff\quad
\mathbf x-\mathbf y\in\Z\mathbf d
\]
because the entries of $\mathbf x-\mathbf y$ sum to $h-h=0$ and
$\mathbf d\cdot\vec{A}=0$. Since $hA=\{\mathbf x\cdot\vec{A}:\mathbf x\in\X_{h,3}\}$,
$\abs{hA}$ is the number of classes of this equivalence relation. The class of
$\mathbf y$ is $\{\mathbf y+t\mathbf d:t\in T_{\mathbf y}\}$, where, because
$d_1+d_2+d_3=0$,
\[
T_{\mathbf y}
 \ = \ \{t\in\Z:\mathbf y+t\mathbf d\in\X_{h,3}\}
 \ = \ \bigcap_{i=1}^3\{t\in\Z:y_i+td_i\geq0\}
 \ = \ \Bigl[\max_{d_i>0}\Bigl\lceil-\frac{y_i}{d_i}\Bigr\rceil
        ..\min_{d_i<0}\Bigl\lfloor\frac{y_i}{\abs{d_i}}\Bigr\rfloor\Bigr].
\]
Both index sets are nonempty because $\mathbf d\neq\mathbf 0$ and
$d_1+d_2+d_3=0$. Thus each class $C$ equals
$\{\mathbf c,\mathbf c+\mathbf d,\ldots,\mathbf c+\alpha_C\mathbf d\}$, so
$\abs C=\alpha_C+1$ and $C$ contains exactly $\alpha_C$ pairs $(\mathbf x,\mathbf y)$
with $\mathbf x-\mathbf y=\mathbf d$. Every such pair lies in a single class, so with
\[
N \ = \ \#\{(\mathbf x,\mathbf y)\in\X_{h,3}\times\X_{h,3}:
\mathbf x-\mathbf y=\mathbf d\}
\]
we get
\begin{align*}
\abs{hA}
& \ = \ \sum_C1
 \ = \ \sum_C\bigl((\alpha_C+1)-\alpha_C\bigr)\\
& \ = \ \sum_C\abs C-\sum_C\alpha_C
 \ = \ \abs{\X_{h,3}}-N
 \ = \ M_{h,3}-N.
\end{align*}

Since $\mathbf d\in\Rcal_{h,3}$, we have $s\leq h$. As in the proof of
Lemma~\ref{lem:equal-sum-criterion}, the map
\[
\phi(\mathbf t) \ = \ (\mathbf d^++\mathbf t,\ \mathbf d^-+\mathbf t)
\]
sends $\X_{h-s,3}$ into the set of pairs counted by $N$. Conversely, for a pair
$(\mathbf x,\mathbf y)$ counted by $N$, let
$\mathbf t=\psi(\mathbf x,\mathbf y)=\mathbf x-\mathbf d^+$. Then
\begin{align*}
\mathbf x-\mathbf y \ = \ \mathbf d \ = \ \mathbf d^+-\mathbf d^-
& \ \implies\ \mathbf t \ = \ \mathbf x-\mathbf d^+ \ = \ \mathbf y-\mathbf d^-,\\
d_i\leq0
& \ \implies\ d_i^+=0
  \ \implies\ t_i \ = \ x_i \ \geq\ 0,\\
d_i>0
& \ \implies\ d_i^-=0
  \ \implies\ t_i \ = \ y_i \ \geq\ 0,\\
\sum_{i=1}^3t_i
& \ = \ \sum_{i=1}^3x_i-\sum_{i=1}^3d_i^+
  \ = \ h-w(\mathbf d) \ = \ h-s,
\end{align*}
so $\mathbf t\in\X_{h-s,3}$. Moreover,
\begin{align*}
\psi(\phi(\mathbf t))
& \ = \ (\mathbf d^++\mathbf t)-\mathbf d^+ \ = \ \mathbf t,\\
\phi(\psi(\mathbf x,\mathbf y))
& \ = \ (\mathbf x,\ \mathbf x-\mathbf d^++\mathbf d^-)
  \ = \ (\mathbf x,\ \mathbf x-\mathbf d) \ = \ (\mathbf x,\mathbf y).
\end{align*}
Hence $\phi$ is a bijection, and
\[
N \ = \ \abs{\X_{h-s,3}} \ = \ M_{h-s,3},
\qquad
\abs{hA} \ = \ M_{h,3}-M_{h-s,3}. \qedhere
\]
\end{proof}

\begin{example}\label{ex:E-size}
By Example~\ref{ex:E-relation}, for every $h\geq2$ the vector $(1,-2,1)$, of weight $2$, is the only vector in $\Rcal_{h,3}$ orthogonal to $\vec{E}$. By Lemma~\ref{lem:gap-size-lemma}, we get
\[
\abs{hE} \ = \ M_{h,3}-M_{h-2,3} \ = \ \binom{h+2}2-\binom h2 \ = \ 2h+1,
\]
in agreement with $hE=[h..3h]$. For $h=4$, this recovers $\abs{4E}=9$.
\end{example}

\section{Main Theorem}

\maintheorem*

\begin{proof}[Proof of Theorem~\ref{thm:main}]
For $\mathbf d\in\Rcal_{h,k}$, let
\[
E_{\mathbf d}
 \ = \ \left\{A\in\binom{[q]}k:\mathbf d\cdot\vec{A}=0\right\}.
\]
By Lemma~\ref{lem:equal-sum-criterion},
\[
\abs{hA}<M_{h,k}
\quad\iff\quad
A\in\bigcup_{\mathbf d\in\Rcal_{h,k}}E_{\mathbf d}.
\]
Using Lemma~\ref{lem:rank-lemma} and
Proposition~\ref{prop:coefficient-vectors-lemma}, we obtain
\begin{equation}\label{eq:not-Bh-count}
\begin{aligned}
\#\left\{A\in\binom{[q]}k:\abs{hA}<M_{h,k}\right\}
& \ = \ \#\bigcup_{\mathbf d\in\Rcal_{h,k}}E_{\mathbf d}\\
& \ \leq\ \sum_{\mathbf d\in\Rcal_{h,k}}\abs{E_{\mathbf d}}\\
& \ \leq\ \abs{\Rcal_{h,k}}q^{k-1}\\
& \ \leq\ (2h+1)^{k-1}q^{k-1}\\
& \ = \ O(h^{k-1}q^{k-1}).
\end{aligned}
\end{equation}
Consequently,
\begin{align*}
\#\left\{A\in\binom{[q]}k:\abs{hA}=M_{h,k}\right\}
& \ = \ \binom qk-
\#\left\{A\in\binom{[q]}k:\abs{hA}<M_{h,k}\right\}\\
& \ = \ \binom qk+O(h^{k-1}q^{k-1})\\
& \ = \ \frac{q^k}{k!}+O(h^{k-1}q^{k-1})\\
& \ = \ \Theta(q^k),
\end{align*}
which proves~(1). The last step uses $h^{k-1}\leq h^k=o(q)$, which holds because
$h=o(q^{1/k})$.

For~(2), let
\[
\eta(A) \ = \ \#\{\mathbf d\in\Rcal_{h,k}:\mathbf d\cdot\vec{A}=0\},
\]
the number of coefficient vectors in $\Rcal_{h,k}$ orthogonal to $\vec{A}$.
Lemma~\ref{lem:equal-sum-criterion} gives
\[
\eta(A)=0
\quad\implies\quad
\abs{hA}=M_{h,k}=M_{h,k}-M_{-1,k}=P_0.
\]
If $\eta(A)=1$, let $\mathbf d$ be the unique vector and let $s=w(\mathbf d)$.
Weight $1$ would give entries $1$ and $-1$ in distinct coordinates $i,j$,
and hence
\[
0=\mathbf d\cdot\vec{A}=a_i-a_j\neq0.
\]
Thus $2\leq s\leq h$, and Lemma~\ref{lem:gap-size-lemma} gives
\[
\abs{hA} \ = \ M_{h,k}-M_{h-s,k} \ = \ P_{h-s+1},
\qquad 1\leq h-s+1\leq h-1.
\]
It follows that
\[
\abs{hA}\notin\{P_0,\ldots,P_{h-1}\}
\quad\implies\quad
\eta(A)\geq2.
\]
For distinct $\mathbf d,\mathbf e\in\Rcal_{h,k}$,
Proposition~\ref{prop:coefficient-vectors-lemma} and
Lemma~\ref{lem:rank-lemma} give
\[
\abs{E_{\mathbf d}\cap E_{\mathbf e}}
\ \leq\
\#\{\mathbf v\in[q]^k:
\mathbf d\cdot\mathbf v=\mathbf e\cdot\mathbf v=0\}
\ \leq\ q^{k-2}.
\]
Therefore
\begin{align*}
\#\left\{A\in\binom{[q]}k:
\abs{hA}\notin\{P_0,\ldots,P_{h-1}\}\right\}
& \ \leq\ \sum_{\substack{\mathbf d,\mathbf e\in\Rcal_{h,k}\\
                         \mathbf d\neq\mathbf e}}
\abs{E_{\mathbf d}\cap E_{\mathbf e}}\\
& \ \leq\ \abs{\Rcal_{h,k}}^2q^{k-2}\\
& \ \leq\ (2h+1)^{2k-2}q^{k-2}\\
& \ = \ O(h^{2k-2}q^{k-2}),
\end{align*}
proving~(2).

For~(3), fix $1\leq\ell<h$ and set $s=h-\ell+1$, so $2\leq s\leq h$.
Choose
\[
1\leq a\leq\left\lfloor\frac q2\right\rfloor,
\qquad
1\leq r\leq\left\lfloor\frac q{2s}\right\rfloor,
\qquad
(a_1,a_2,a_3)=(a,a+r,a+sr).
\]
Then
\[
1\leq a_1<a_2<a_3
\leq\left\lfloor\frac q2\right\rfloor
   +s\left\lfloor\frac q{2s}\right\rfloor
\leq q.
\]
Let $\mathbf d_0=(s-1,-s,1,0,\ldots,0)$. We have
\begin{align*}
\sum_{i=1}^k(d_0)_i
& \ = \ (s-1)-s+1 \ = \ 0,\\
\gcd(\abs{(d_0)_1},\ldots,\abs{(d_0)_k})
& \ = \ \gcd(s-1,s,1) \ = \ 1,\\
w(\mathbf d_0)
& \ = \ (s-1)+1 \ = \ s \ \leq\ h,\\
(d_0)_{j(\mathbf d_0)}
& \ = \ s-1 \ > \ 0.
\end{align*}
Thus $\mathbf d_0\in\Rcal_{h,k}$. For every extension to a labeled tuple
$\mathbf a=(a_1,\ldots,a_k)$,
\begin{align*}
\mathbf d_0\cdot\mathbf a
& \ = \ (s-1)a_1-sa_2+a_3\\
& \ = \ (s-1)a-s(a+r)+(a+sr)\\
& \ = \ ((s-1)-s+1)a+(-s+s)r \ = \ 0.
\end{align*}

Suppose that $\mathbf e=(e_1,e_2,e_3,0,\ldots,0)\in\Rcal_{h,k}$ and
$\mathbf e\cdot\mathbf a=0$. Then
\begin{align*}
0
& \ = \ e_1a+e_2(a+r)+e_3(a+sr)\\
& \ = \ a(e_1+e_2+e_3)+r(e_2+se_3)\\
& \ = \ r(e_2+se_3).
\end{align*}
Since $r>0$ and $e_1+e_2+e_3=0$,
\[
e_2=-se_3,
\qquad
e_1=(s-1)e_3,
\qquad
\mathbf e=e_3\mathbf d_0.
\]
Primitivity and the sign convention give
\[
1=\gcd(\abs{e_1},\ldots,\abs{e_k})=\abs{e_3},
\qquad
0<e_1=(s-1)e_3
\quad\implies\quad
e_3=1,
\]
so $\mathbf e=\mathbf d_0$.

Suppose that $a_1,\ldots,a_{j-1}$ are distinct and that $\mathbf d_0$ is the
unique vector supported on the first $j-1$ coordinates with zero dot
product. For $\mathbf e\in\Rcal_{h,k}$ supported on the first $j$
coordinates and satisfying $e_j\neq0$,
\[
\sum_{i=1}^j e_i a_i=0
\quad\iff\quad
a_j=-\frac1{e_j}\sum_{i=1}^{j-1}e_i a_i.
\]
Each such vector excludes at most one value of $a_j$, because $e_j\neq0$ and the
equivalence above determines $a_j$ from $a_1,\ldots,a_{j-1}$. Since every $e_i$
lies in $[-h,h]$, at most $(2h+1)^j$ values are excluded by these vectors.
If $F_j$ is the union of these values and $a_1,\ldots,a_{j-1}$, then
\begin{align*}
\#\{a_j\in[q]:a_j\notin F_j\}
& \ = \ q-\abs{F_j\cap[q]}\\
& \ \geq\ q-(2h+1)^j-(j-1)\\
& \ \geq\ \frac q2
\end{align*}
for all sufficiently large $q$, because $(2h+1)^j+j\leq(3h)^k+k$ and $h^k=o(q)$. Vectors with $e_j=0$ are covered by the
inductive assumption, and vectors with $e_j\neq0$ are excluded by the
choice of $a_j$. Repeating for $j=4,\ldots,k$ therefore gives
\[
\{\mathbf e\in\Rcal_{h,k}:\mathbf e\cdot\mathbf a=0\}
 \ = \ \{\mathbf d_0\}.
\]
The number of labeled tuples produced is at least
\begin{align*}
\left\lfloor\frac q2\right\rfloor
\left\lfloor\frac q{2s}\right\rfloor
\prod_{j=4}^k\frac q2
& \ = \ \left\lfloor\frac q2\right\rfloor
\left\lfloor\frac q{2s}\right\rfloor
\left(\frac q2\right)^{k-3}\\
& \ = \ \Omega\!\left(\frac{q^{k-1}}s\right).
\end{align*}
Here we used $\lfloor x\rfloor\geq x/2$ for $x\geq1$, which applies to
$x=q/(2s)$ because $s\leq h=o(q)$.

Let $A$ be the set of entries of $\mathbf a$, and let $P$ be the permutation
matrix with $\vec{A}=P\mathbf a$. Choose the prescribed sign
$\varepsilon$ and let $\mathbf d'=\varepsilon P\mathbf d_0$. Then
\begin{align*}
\mathbf d'\cdot\vec{A}
& \ = \ (\varepsilon P\mathbf d_0)\cdot(P\mathbf a)
  \ = \ \varepsilon\mathbf d_0\cdot\mathbf a \ = \ 0,\\
w(\mathbf d')
& \ = \ w(\mathbf d_0) \ = \ s,
\qquad \mathbf d'\in\Rcal_{h,k}.
\end{align*}
For any $\mathbf e\in\Rcal_{h,k}$ with $\mathbf e\cdot\vec{A}=0$, choose
$\delta\in\{1,-1\}$ so that $\delta P^{-1}\mathbf e\in\Rcal_{h,k}$. Then
\[
(\delta P^{-1}\mathbf e)\cdot\mathbf a
 \ = \ \delta\mathbf e\cdot\vec{A} \ = \ 0
\quad\implies\quad
\delta P^{-1}\mathbf e=\mathbf d_0
\quad\implies\quad
\mathbf e=\delta P\mathbf d_0.
\]
The sign convention forces $\mathbf e=\mathbf d'$. Hence
Lemma~\ref{lem:gap-size-lemma} gives
\begin{align*}
\abs{hA}
& \ = \ M_{h,k}-M_{h-s,k}\\
& \ = \ M_{h,k}-M_{\ell-1,k}\\
& \ = \ P_\ell.
\end{align*}
At most $k!$ labeled tuples produce the same set, so
\begin{align*}
\#\left\{A\in\binom{[q]}k:\abs{hA}=P_\ell\right\}
& \ \geq\ \frac{\#\{\text{labeled tuples produced}\}}{k!}\\
& \ = \ \Omega\!\left(\frac{q^{k-1}}s\right)
  \ = \ \Omega\!\left(\frac{q^{k-1}}{h-\ell+1}\right).
\end{align*}
Finally, $P_0=M_{h,k}$ and~(1) give
\[
\#\left\{A\in\binom{[q]}k:\abs{hA}=P_0\right\}
 \ = \ \Theta(q^k)
 \ = \ \Omega\!\left(\frac{q^{k-1}}{h+1}\right).
\]
This proves~(3).
\end{proof}

\begin{corollary}\label{cor:gap-pattern}
Fix $k\geq3$, and let $h=h(q)\geq2$ satisfy $h=o(q^{1/(2k-1)})$. For all
sufficiently large $q$, the $h$ most
frequent values of $\abs{hA}$, with $A\in\binom{[q]}k$, are precisely
\[
P_0,P_1,\ldots,P_{h-1}.
\]
Moreover,
\begin{equation}
P_\ell-P_{\ell+1}
 \ = \ \binom{\ell+k-2}{k-2},
\qquad 0\leq\ell<h-1.
\end{equation}
In particular, the gaps are consecutive triangular numbers when $k=4$.
\end{corollary}

\begin{proof}
For an integer $v$, let
\[
f_q(v) \ = \ \#\left\{A\in\binom{[q]}k:\abs{hA}=v\right\}.
\]
Since $h=o(q^{1/(2k-1)})$ implies $h=o(q^{1/k})$, Theorem~\ref{thm:main} applies and
gives
\begin{align*}
f_q(P_0)& \ = \ \Theta(q^k),\\
f_q(P_\ell)& \ = \ \Omega\!\left(\frac{q^{k-1}}{h-\ell+1}\right)
 \ = \ \Omega\!\left(\frac{q^{k-1}}h\right)
 \qquad(1\leq\ell<h),\\
\sum_{v\notin\{P_0,\ldots,P_{h-1}\}}f_q(v)
& \ = \ O(h^{2k-2}q^{k-2}).
\end{align*}
Thus there are constants $c,C>0$ such that, for all sufficiently large $q$,
\begin{align*}
0\leq
\frac{\max_{v\notin\{P_0,\ldots,P_{h-1}\}}f_q(v)}
     {\min_{0\leq\ell<h}f_q(P_\ell)}
& \ \leq\
\frac{\sum_{v\notin\{P_0,\ldots,P_{h-1}\}}f_q(v)}
     {\min_{0\leq\ell<h}f_q(P_\ell)}\\
& \ \leq\ \frac{Ch^{2k-2}q^{k-2}}{cq^{k-1}/h}\\
& \ = \ \frac Cc\frac{h^{2k-1}}q
 \ \longrightarrow\ 0,
\end{align*}
because $h^{2k-1}=o(q)$.
Consequently, each listed value occurs more frequently than every
unlisted value for all sufficiently large $q$.

For $0\leq\ell<h-1$, Pascal's identity gives
\begin{align*}
P_\ell-P_{\ell+1}
& \ = \ (M_{h,k}-M_{\ell-1,k})-(M_{h,k}-M_{\ell,k})\\
& \ = \ M_{\ell,k}-M_{\ell-1,k}\\
& \ = \ \binom{\ell+k-1}{k-1}-\binom{\ell+k-2}{k-1}\\
& \ = \ \binom{\ell+k-2}{k-2}.
\end{align*}
For $\ell=0$, this uses $M_{-1,k}=0$ and $M_{0,k}=1$. For $k=4$, the gaps are
\[
\binom{\ell+2}{2} \ = \ \frac{(\ell+1)(\ell+2)}2,
\]
the consecutive triangular numbers.
\end{proof}

\begin{remark}\label{rem:explicit}
The proofs of Theorem~\ref{thm:main} and Corollary~\ref{cor:gap-pattern} give an explicit value of $q$ beyond which the conclusion of Corollary~\ref{cor:gap-pattern} holds. Every set counted in Theorem~\ref{thm:main}(2) satisfies $\mathbf d\cdot\vec{A}=\mathbf e\cdot\vec{A}=0$ for two distinct $\mathbf d,\mathbf e\in\Rcal_{h,k}$, so by Lemma~\ref{lem:rank-lemma} and Proposition~\ref{prop:coefficient-vectors-lemma}, there are at most $(2h+1)^{2k-2}q^{k-2}/2$ such sets for every $q$. In the same way, at most $(2h+1)^{k-1}q^{k-1}$ sets satisfy $\abs{hA}<M_{h,k}$. If $q\geq2((2h+1)^k+k-1)$, then the construction in the proof of Theorem~\ref{thm:main}(3) produces at least
\[
\frac{1}{k!}\left\lfloor\frac q2\right\rfloor\left\lfloor\frac q{2s}\right\rfloor\left(\frac q2\right)^{k-3}\geq\frac{q^{k-1}}{2^{k+1}s\,k!}
\]
sets with $\abs{hA}=P_\ell$, where $s=h-\ell+1\leq h$. Comparing these bounds, and using $\binom qk\geq q^k/(2^kk!)$ for $q\geq2k$ in the case $\ell=0$, we find that the conclusion of Corollary~\ref{cor:gap-pattern} holds whenever
\[
q>2^k\,k!\,h\,(2h+1)^{2k-2}.
\]
In particular, $h$ may grow with $q$, provided that $2^k\,k!\,3^{2k-2}h^{2k-1}<q$. This bound is far from sharp. For example, a direct computation shows that for $(h,k)=(3,4)$ the conclusion of Corollary~\ref{cor:gap-pattern} holds for every $q$ with $16\leq q\leq50$, whereas the bound above requires $q>2^4\cdot4!\cdot3\cdot7^6\approx1.4\times10^8$.
\end{remark}

\section{Additional Results and Future Work}\label{sec:additional}

Throughout this section, $h\geq2$ is fixed, so the constants $C$, $c$ and $q_0$
may also depend on $h$.

\subsection{Generalizations of Senger's Results}
This section collects generalizations to arbitrary $k\geq3$ of several results of Senger~\cite{Senger}, though none are needed for the proof of Theorem~\ref{thm:main}. Let
\[
\Bcal_{h,k}(q) \ = \ \left\{A\in\binom{[q]}k:A\text{ is a }B_h\text{-set}\right\}
\]
and
\[
\Bcal_{h,k}^*(q) \ = \ \left\{A\in\binom{[q]}k:A\text{ is a }B_h\text{-set but not a }B_{h+1}\text{-set}\right\}.
\]

\begin{corollary}[{cf.~\cite[Thm.~1]{Senger}}]\label{cor:proportion}
Fix $k\geq3$ and $h\geq2$. Then, as $q\to\infty$,
\[
1-\frac{\abs{\Bcal_{h,k}(q)}}{\binom qk}
 \ = \ \Theta(q^{-1}).
\]
In particular,
\[
\frac{\abs{\Bcal_{h,k}(q)}}{\binom qk}\longrightarrow1.
\]
For $(h,k)=(2,3)$,
\[
\lim_{q\to\infty}q\left(1-
\frac{\abs{\Bcal_{2,3}(q)}}{\binom q3}\right) \ = \ \frac32.
\]
\end{corollary}

\begin{proof}
Let $N_q=\binom qk-\abs{\Bcal_{h,k}(q)}$. By the definition of a $B_h$-set,
\begin{align*}
N_q
& \ = \ \#\left\{A\in\binom{[q]}k:\abs{hA}<M_{h,k}\right\},\\
P_1
& \ = \ M_{h,k}-M_{0,k} \ = \ M_{h,k}-1.
\end{align*}
Theorem~\ref{thm:main}(3) and \eqref{eq:not-Bh-count} therefore give a
constant $c>0$ such that, for all sufficiently large $q$,
\[
cq^{k-1}
 \ \leq\ \#\left\{A\in\binom{[q]}k:\abs{hA}=P_1\right\}
 \ \leq\ N_q
 \ \leq\ (2h+1)^{k-1}q^{k-1}.
\]
Thus $N_q=\Theta(q^{k-1})$. Since $\binom qk=\Theta(q^k)$,
\begin{align*}
1-\frac{\abs{\Bcal_{h,k}(q)}}{\binom qk}
& \ = \ \frac{\binom qk-\abs{\Bcal_{h,k}(q)}}{\binom qk}\\
& \ = \ \frac{N_q}{\binom qk}
 \ = \ \Theta(q^{-1})
 \ \longrightarrow\ 0.
\end{align*}

For $(h,k)=(2,3)$ and $a_1<a_2<a_3$, the six possible two-fold sums satisfy
\[
2a_1<a_1+a_2<\min\{2a_2,a_1+a_3\}
\leq\max\{2a_2,a_1+a_3\}<a_2+a_3<2a_3.
\]
Hence the only possible collision is $2a_2=a_1+a_3$, and
\[
\binom{[q]}3\setminus\Bcal_{2,3}(q)
 \ = \ \{\{a,a+r,a+2r\}:a,r\geq1,\ a+2r\leq q\}.
\]
Let $m=\lfloor(q-1)/2\rfloor$. Then
\begin{align*}
\binom q3-\abs{\Bcal_{2,3}(q)}
& \ = \ \sum_{r=1}^m(q-2r)\\
& \ = \ mq-m(m+1)\\
& \ = \ \frac{q^2}{4}+O(q),\\
q\left(1-\frac{\abs{\Bcal_{2,3}(q)}}{\binom q3}\right)
& \ = \ \frac{q(q^2/4+O(q))}{q(q-1)(q-2)/6}\\
& \ = \ \frac{3/2+O(q^{-1})}{(1-q^{-1})(1-2q^{-1})}
 \ \longrightarrow\ \frac32.
\end{align*}
\end{proof}

For $h=2$, the convergence to $1$ also follows from work of
Nathanson~\cite{Nathanson2004}, which allows $k$ to grow with $q$.

\begin{remark}\label{rem:k3-limit}
For $k=3$, the limit in Corollary~\ref{cor:proportion} can be computed for every $h\geq2$. Let $A\in\binom{[q]}3$ and write
\[
\frac{a_2-a_1}{a_3-a_1} \ = \ \frac uw,
\qquad \gcd(u,w)=1,
\qquad 0<u<w.
\]
If $\mathbf d\in\Z^3\setminus\{\mathbf 0\}$ has $d_1+d_2+d_3=0$, then
\[
\mathbf d\cdot\vec{A} \ = \ d_2(a_2-a_1)+d_3(a_3-a_1),
\]
so $\mathbf d\cdot\vec{A}=0$ if and only if $d_2u=-d_3w$, that is, $\mathbf d$ is an integer multiple of $(w-u,-w,u)$. This vector is primitive, has weight $(w-u)+u=w$, and has positive first coordinate. By Lemma~\ref{lem:equal-sum-criterion}, $A$ is a $B_h$-set if and only if $w>h$. For fixed coprime $0<u<w$, the sets with $(a_2-a_1)/(a_3-a_1)=u/w$ are $\{a,a+ut,a+wt\}$ with $a,t\geq1$ and $a+wt\leq q$, and there are
\[
\sum_{t=1}^{\lfloor(q-1)/w\rfloor}(q-wt) \ = \ \frac{q^2}{2w}+O(q)
\]
of them. Summing over the $\varphi(w)$ admissible values of $u$ for each $2\leq w\leq h$, where $\varphi$ is Euler's totient function, gives
\[
\binom q3-\abs{\Bcal_{h,3}(q)} \ = \ \frac{q^2}2\sum_{w=2}^h\frac{\varphi(w)}w+O(h^2q),
\]
and therefore
\[
\lim_{q\to\infty}q\left(1-\frac{\abs{\Bcal_{h,3}(q)}}{\binom q3}\right) \ = \ 3\sum_{w=2}^h\frac{\varphi(w)}w.
\]
For $h=2$ this recovers the value $3/2$ in Corollary~\ref{cor:proportion}, and for $h=3$ it gives $7/2$.
\end{remark}

\begin{lemma}[{cf.~\cite[Lemma~4]{Senger}}]\label{lem:disjoint-support}
Given $A\in\Bcal_{h,k}(q)$, any distinct sum vectors
$\mathbf x,\mathbf y\in\X_{h+1,k}$ that have the same dot product with
$\vec{A}$ must have disjoint support, and therefore are orthogonal.
\end{lemma}

Senger proved Lemma~\ref{lem:disjoint-support} for $k=4$~\cite[Lemma~4]{Senger}. His proof uses only the fact that $A$ is a $B_h$-set, so it applies without change for every $k\geq3$.

If $A\in\Bcal_{h,k}^*(q)$, then $A$ is not a $B_{h+1}$-set, so there are distinct $\mathbf x,\mathbf y\in\X_{h+1,k}$ with $\mathbf x\cdot\vec{A}=\mathbf y\cdot\vec{A}$. By Lemma~\ref{lem:disjoint-support}, every such pair has disjoint support, because $A$ is a $B_h$-set. For distinct $\mathbf x,\mathbf y\in\X_{h+1,k}$ with disjoint support, let
\[
T_{\mathbf x,\mathbf y}
 \ = \ \#\{A\in\Bcal_{h,k}^*(q):
\mathbf x\cdot\vec{A}=\mathbf y\cdot\vec{A}\}.
\]
Then
\[
\Bcal_{h,k}^*(q) \ = \ \bigcup_{(\mathbf x,\mathbf y)}
\left\{A\in\Bcal_{h,k}^*(q):\mathbf x\cdot\vec{A}=\mathbf y\cdot\vec{A}\right\},
\]
where the union runs over all pairs of distinct $\mathbf x,\mathbf y\in\X_{h+1,k}$ with disjoint support. The following proposition determines when each $T_{\mathbf x,\mathbf y}$ is zero and when it has order $q^{k-1}$.

\begin{proposition}\label{prop:T-dichotomy}
Let $\mathbf x,\mathbf y\in\X_{h+1,k}$ be distinct sum vectors with disjoint support, let $\mathbf d=\mathbf x-\mathbf y$, and let $D_j=d_j+d_{j+1}+\cdots+d_k$ for $2\leq j\leq k$. Then $T_{\mathbf x,\mathbf y}\leq q^{k-1}$ for every $q$. Moreover, the following statements hold.
\begin{enumerate}
\item If $\mathbf d$ is not primitive, or if $D_2,\ldots,D_k$ are all nonnegative or all nonpositive, then $T_{\mathbf x,\mathbf y}=0$ for every $q$.
\item Otherwise, $T_{\mathbf x,\mathbf y}=\Theta(q^{k-1})$.
\end{enumerate}
\end{proposition}

\begin{proof}
Disjoint support gives
\[
\mathbf d^+=\mathbf x,
\qquad \mathbf d^-=\mathbf y,
\qquad
w(\mathbf d)=\sum_{i=1}^k x_i=h+1,
\qquad
\sum_{i=1}^k d_i=(h+1)-(h+1)=0.
\]
Since $\mathbf d\neq\mathbf 0$, Lemma~\ref{lem:rank-lemma} gives
\begin{align*}
T_{\mathbf x,\mathbf y}
& \ = \ \#\{A\in\Bcal_{h,k}^*(q):\mathbf d\cdot\vec{A}=0\}\\
& \ \leq\ \#\{\mathbf v\in[q]^k:\mathbf d\cdot\mathbf v=0\}\\
& \ \leq\ q^{k-1}.
\end{align*}

For~(1), first suppose that $g=\gcd(\abs{d_1},\ldots,\abs{d_k})\geq2$,
and let $\mathbf d'=\mathbf d/g$. Then
\begin{align*}
\gcd(\abs{d'_1},\ldots,\abs{d'_k})
& \ = \ \frac gg \ = \ 1,\\
\sum_{i=1}^k d'_i
& \ = \ \frac1g\sum_{i=1}^k d_i \ = \ 0,\\
w(\mathbf d')
& \ = \ \frac{h+1}{g}\leq\frac{h+1}{2}\leq h.
\end{align*}
Thus $\mathbf d'$ or $-\mathbf d'$ lies in $\Rcal_{h,k}$. By
Lemma~\ref{lem:equal-sum-criterion},
\[
\mathbf d\cdot\vec{A}=0
\quad\implies\quad
\mathbf d'\cdot\vec{A}=\frac1g\mathbf d\cdot\vec{A}=0
\quad\implies\quad
\abs{hA}<M_{h,k}.
\]
Therefore
\[
\{A\in\Bcal_{h,k}^*(q):\mathbf d\cdot\vec{A}=0\}
 \ = \ \varnothing,
\qquad T_{\mathbf x,\mathbf y}=0.
\]

Next, write $t_1=a_1$ and $t_j=a_j-a_{j-1}$ for $2\leq j\leq k$. Then
$t_j>0$ and $a_i=\sum_{j=1}^i t_j$, so
\begin{equation}\label{eq:suffix-sums}
\begin{aligned}
\mathbf d\cdot\vec{A}
& \ = \ \sum_{i=1}^k d_i\sum_{j=1}^i t_j\\
& \ = \ \sum_{j=1}^k t_j\sum_{i=j}^k d_i\\
& \ = \ t_1\sum_{i=1}^k d_i+\sum_{j=2}^k t_jD_j\\
& \ = \ \sum_{j=2}^k t_jD_j.
\end{aligned}
\end{equation}
The identities
\[
d_1=-D_2,
\qquad d_j=D_j-D_{j+1}\quad(2\leq j<k),
\qquad d_k=D_k
\]
give
\[
D_2=\cdots=D_k=0\quad\implies\quad\mathbf d=\mathbf 0,
\]
which is impossible. Thus \eqref{eq:suffix-sums} gives
\begin{align*}
D_j\geq0\text{ for every }j
& \ \implies\ \mathbf d\cdot\vec{A}=\sum_{j=2}^k t_jD_j>0,\\
D_j\leq0\text{ for every }j
& \ \implies\ \mathbf d\cdot\vec{A}=\sum_{j=2}^k t_jD_j<0.
\end{align*}
In either case $\mathbf d\cdot\vec{A}\neq0$ for every $A$, and hence
$T_{\mathbf x,\mathbf y}=0$.

For~(2), suppose that $\mathbf d$ is primitive and that the $D_j$ take
both signs. Let
\[
L_+=\sum_{\substack{2\leq j\leq k\\D_j>0}}D_j>0,
\qquad
L_-=-\sum_{\substack{2\leq j\leq k\\D_j<0}}D_j>0.
\]
Choose $\widetilde\tau_1=1$ and, for $j\geq2$, set
\[
\widetilde\tau_j=
\begin{cases}
L_- & \text{if }D_j>0,\\
L_+ & \text{if }D_j<0,\\
1 & \text{if }D_j=0.
\end{cases}
\]
Then
\[
\sum_{j=2}^k\widetilde\tau_jD_j
 \ = \ L_-L_+-L_+L_- \ = \ 0.
\]
Let $\tau_j=\widetilde\tau_j/(2\sum_{i=1}^k\widetilde\tau_i)$ and
$u_i=\sum_{j=1}^i\tau_j$. We have
\begin{align*}
u_1& \ = \ \tau_1>0,\\
u_{i+1}-u_i& \ = \ \tau_{i+1}>0,\\
u_k& \ = \ \sum_{j=1}^k\tau_j \ = \ \frac12<1,\\
\mathbf d\cdot\mathbf u
& \ = \ \sum_{j=2}^k\tau_jD_j
  \ = \ \frac{\sum_{j=2}^k\widetilde\tau_jD_j}
                 {2\sum_{i=1}^k\widetilde\tau_i} \ = \ 0.
\end{align*}
Hence $\mathbf u\in U\cap H$, where
\[
U=\{\mathbf v\in\mathbb R^k:0<v_1<\cdots<v_k<1\},
\qquad
H=\{\mathbf v\in\mathbb R^k:\mathbf d\cdot\mathbf v=0\}.
\]
Since $U$ is open, choose $\rho>0$ so that the ball
$B=\{\mathbf v\in H:\lVert\mathbf v-\mathbf u\rVert<\rho\}$ lies in
$U\cap H$.

The integer kernel $\Lambda=\Z^k\cap H$ is a lattice. Primitivity gives
\[
\mathbf d\cdot\Z^k=\Z,
\qquad
\operatorname{rank}(\Lambda)
 \ = \ k-\operatorname{rank}(\Z) \ = \ k-1 \ = \ \dim H.
\]
Choose a lattice basis $\mathbf b_1,\ldots,\mathbf b_{k-1}$ and let
$L=\sum_{i=1}^{k-1}\lVert\mathbf b_i\rVert$. Write
$q\mathbf u=\sum_{i=1}^{k-1}\theta_i\mathbf b_i$, and let
\[
\boldsymbol\lambda_q
 \ = \ \sum_{i=1}^{k-1}\left\lfloor\theta_i+\frac12\right\rfloor\mathbf b_i.
\]
Then
\[
\lVert\boldsymbol\lambda_q-q\mathbf u\rVert
\leq\sum_{i=1}^{k-1}
\left\lvert\left\lfloor\theta_i+\frac12\right\rfloor-\theta_i\right\rvert
\lVert\mathbf b_i\rVert
\leq\frac L2.
\]
Let $m=\lfloor q\rho/(2L)\rfloor$. For integers $z_i$ with $\abs{z_i}\leq m$,
\[
\left\lVert\boldsymbol\lambda_q+\sum_{i=1}^{k-1}z_i\mathbf b_i
             -q\mathbf u\right\rVert
\leq\frac L2+mL
\leq\frac L2+\frac{q\rho}{2}
<q\rho
\]
for all sufficiently large $q$. All these distinct lattice points lie in
$qB$, so there is $c>0$ such that
\[
\#(\Lambda\cap qB)
\geq(2m+1)^{k-1}
=\left(2\left\lfloor\frac{q\rho}{2L}\right\rfloor+1\right)^{k-1}
\geq cq^{k-1}.
\]
Each $\mathbf v\in\Lambda\cap qB$ satisfies
\[
0<v_1<\cdots<v_k<q,
\qquad
\mathbf d\cdot\mathbf v=0,
\]
so $\mathbf v=\vec{A}$ for a set $A\in\binom{[q]}k$, and
\[
\mathbf x\cdot\vec{A}-\mathbf y\cdot\vec{A}
 \ = \ \mathbf d\cdot\vec{A} \ = \ 0.
\]
Since $\mathbf x\neq\mathbf y$, each such $A$ is not a $B_{h+1}$-set.

It remains to discard those that are not $B_h$-sets. For
$\mathbf e\in\Rcal_{h,k}$,
\[
w(\mathbf e)\leq h<h+1=w(\mathbf d)
\quad\implies\quad\mathbf e\neq\pm\mathbf d.
\]
The primitive-vector argument in
Proposition~\ref{prop:coefficient-vectors-lemma} therefore gives linear
independence of $\mathbf d$ and $\mathbf e$. Lemma~\ref{lem:rank-lemma}
and Proposition~\ref{prop:coefficient-vectors-lemma} yield
\begin{align*}
\#\{\text{sets discarded}\}
& \ \leq\ \sum_{\mathbf e\in\Rcal_{h,k}}
\#\{\mathbf v\in[q]^k:
\mathbf d\cdot\mathbf v=\mathbf e\cdot\mathbf v=0\}\\
& \ \leq\ \abs{\Rcal_{h,k}}q^{k-2}\\
& \ \leq\ (2h+1)^{k-1}q^{k-2}.
\end{align*}
Consequently, for all sufficiently large $q$,
\begin{align*}
T_{\mathbf x,\mathbf y}
& \ \geq\ \#(\Lambda\cap qB)-\#\{\text{sets discarded}\}\\
& \ \geq\ cq^{k-1}-(2h+1)^{k-1}q^{k-2}\\
& \ = \ q^{k-1}\left(c-\frac{(2h+1)^{k-1}}q\right)\\
& \ \geq\ \frac c2q^{k-1}.
\end{align*}
Together with $T_{\mathbf x,\mathbf y}\leq q^{k-1}$, this gives
$T_{\mathbf x,\mathbf y}=\Theta(q^{k-1})$.
\end{proof}

\begin{example}\label{ex:T}
Let $\mathbf x=(h,0,1,0,\ldots,0)$ and $\mathbf y=(0,h+1,0,\ldots,0)$. Then $\mathbf d=(h,-h-1,1,0,\ldots,0)$ is primitive, $D_2=-h<0$ and $D_3=1>0$, so $T_{\mathbf x,\mathbf y}=\Theta(q^{k-1})$ by Proposition~\ref{prop:T-dichotomy}. This vector has the same form as the vector $\mathbf d_0$ in the proof of Theorem~\ref{thm:main}, with $s=h+1$. On the other hand, if $\mathbf x=(0,h,1,0,\ldots,0)$ and $\mathbf y=(h+1,0,\ldots,0)$, then $\mathbf d=(-h-1,h,1,0,\ldots,0)$, and $D_j\geq0$ for every $j\geq2$, so $T_{\mathbf x,\mathbf y}=0$ for every $q$. Directly, $\mathbf d\cdot\vec{A}=h(a_2-a_1)+(a_3-a_1)>0$ for every $A$.
\end{example}

\begin{lemma}[{Generalization of~\cite[Lemma~1]{Senger}}]\label{lem:Bstar-count}
Fix $h\geq2$ and $k\geq3$. As $q\to\infty$, the following statements hold.
\begin{enumerate}
\item
\[
\abs{\Bcal_{h,k}^*(q)} \ = \ O(h^{k-1}q^{k-1}).
\]
\item
\[
\abs{\Bcal_{h,k}^*(q)} \ = \ \Omega(q^{k-1}).
\]
\item
\[
\#\{A\in\Bcal_{h,k}^*(q):
\abs{(h+1)A}\leq M_{h+1,k}-2\}
 \ = \ O(h^{2k-2}q^{k-2}).
\]
\end{enumerate}
In particular,
\[
\abs{\Bcal_{h,k}^*(q)} \ = \ \Theta(q^{k-1}),
\qquad
\#\{A\in\Bcal_{h,k}^*(q):
\abs{(h+1)A}=M_{h+1,k}-1\}
 \ = \ \Theta(q^{k-1}).
\]
\end{lemma}

\begin{proof}
For~(1), every $A\in\Bcal_{h,k}^*(q)$ is not a $B_{h+1}$-set. The bound
in \eqref{eq:not-Bh-count}, with $h+1$ in place of $h$, gives
\begin{align}
\abs{\Bcal_{h,k}^*(q)}
& \ \leq\ \#\left\{A\in\binom{[q]}k:
\abs{(h+1)A}<M_{h+1,k}\right\}\notag\\
& \ \leq\ (2h+3)^{k-1}q^{k-1}\label{eq:Bstar-upper}\\
& \ = \ O(h^{k-1}q^{k-1}).\notag
\end{align}

For~(2), use the first pair $\mathbf x,\mathbf y$ in
Example~\ref{ex:T}. Proposition~\ref{prop:T-dichotomy} gives
\[
\abs{\Bcal_{h,k}^*(q)}
\ \geq\ T_{\mathbf x,\mathbf y}
\ =\ \Theta(q^{k-1}),
\]
proving~(2).

For~(3), let
\[
\mathcal N(A)
 \ = \ \{\mathbf d\in\Rcal_{h+1,k}:\mathbf d\cdot\vec{A}=0\}.
\]
By Lemma~\ref{lem:equal-sum-criterion},
\[
A\in\Bcal_{h,k}^*(q)
\quad\implies\quad
\varnothing\neq\mathcal N(A)
\subseteq\Rcal_{h+1,k}\setminus\Rcal_{h,k}
\quad\implies\quad
w(\mathbf d)=h+1\text{ for every }\mathbf d\in\mathcal N(A).
\]
If $\mathcal N(A)=\{\mathbf d\}$, then
Lemma~\ref{lem:gap-size-lemma} gives
\[
\abs{(h+1)A}
 \ = \ M_{h+1,k}-M_{(h+1)-w(\mathbf d),k}
 \ = \ M_{h+1,k}-M_{0,k}
 \ = \ M_{h+1,k}-1.
\]
Consequently,
\[
\abs{(h+1)A}\leq M_{h+1,k}-2
\quad\implies\quad
\abs{\mathcal N(A)}\geq2.
\]
Any two distinct vectors in $\mathcal N(A)$ are linearly independent by
Proposition~\ref{prop:coefficient-vectors-lemma}. Let $F_q$ be the number
of sets in~(3). Lemma~\ref{lem:rank-lemma} gives
\begin{align}
F_q
& \ \leq\ \sum_{\{\mathbf d,\mathbf e\}\in\binom{\Rcal_{h+1,k}}2}
\#\left\{A\in\binom{[q]}k:
\mathbf d\cdot\vec{A}=\mathbf e\cdot\vec{A}=0\right\}\notag\\
& \ \leq\ \binom{\abs{\Rcal_{h+1,k}}}{2}q^{k-2}\notag\\
& \ \leq\ \frac{(2h+3)^{2k-2}}2q^{k-2}\label{eq:Fq-bound}\\
& \ = \ O(h^{2k-2}q^{k-2}).\notag
\end{align}
This proves~(3). For fixed $h$ and $k$, statements~(1) and~(2) give
$\abs{\Bcal_{h,k}^*(q)}=\Theta(q^{k-1})$. Choose $c>0$ so that
$\abs{\Bcal_{h,k}^*(q)}\geq cq^{k-1}$ for all sufficiently large $q$.
Every set in this family has $\abs{(h+1)A}\leq M_{h+1,k}-1$, so by
\eqref{eq:Fq-bound} we get
\begin{align*}
\#\{A\in\Bcal_{h,k}^*(q):\abs{(h+1)A}=M_{h+1,k}-1\}
& \ = \ \abs{\Bcal_{h,k}^*(q)}-F_q\\
& \ \geq\ cq^{k-1}-\frac{(2h+3)^{2k-2}}2q^{k-2}\\
& \ = \ q^{k-1}\left(c-\frac{(2h+3)^{2k-2}}{2q}\right)\\
& \ \geq\ \frac c2q^{k-1}
\end{align*}
for all sufficiently large $q$. The same count is at most
$\abs{\Bcal_{h,k}^*(q)}$, which is at most $(2h+3)^{k-1}q^{k-1}$ by
\eqref{eq:Bstar-upper}, proving that it is $\Theta(q^{k-1})$.
\end{proof}

\begin{proposition}[{cf.~\cite[Lemma~6]{Senger}}]\label{prop:S-interval}
Let $k\geq3$, $h\geq2$ and $q\geq k$, and define
\[
\mathcal S_{h+1,k}(q)
 \ = \ \left\{\mathbf x\cdot\vec{A}:
\mathbf x\in\X_{h+1,k},\ A\in\binom{[q]}k\right\}.
\]
Then
\[
\mathcal S_{h+1,k}(q) \ = \ [(h+1)..(h+1)q],
\qquad
\abs{\mathcal S_{h+1,k}(q)} \ = \ (h+1)q-h.
\]
For an integer $s$, let
\[
\nu_q(s)
 \ = \ \#\left\{(\mathbf x,A)\in\X_{h+1,k}\times\binom{[q]}k:
\mathbf x\cdot\vec{A}=s\right\}.
\]
Then
\[
\nu_q(s) \ = \ \nu_q\bigl((h+1)(q+1)-s\bigr)
\]
for every $s\in\Z$.
\end{proposition}

\begin{proof}
For $\mathbf x\in\X_{h+1,k}$ and $A\in\binom{[q]}k$,
\[
h+1 \ = \ \sum_{i=1}^k x_i
\ \leq\ \sum_{i=1}^k x_i a_i
\ \leq\ q\sum_{i=1}^k x_i
\ = \ (h+1)q.
\]
Thus $\mathcal S_{h+1,k}(q)\subseteq[(h+1)..(h+1)q]$.

Conversely, let $n\in[(h+1)..(h+1)q]$. By division with remainder,
\[
n=(h+1)a+b,
\qquad 1\leq a\leq q,
\qquad 0\leq b\leq h.
\]
Write $\mathbf e_i$ for the $i$th coordinate vector. If $b=0$, choose
$A\in\binom{[q]}k$ containing $a=a_i$, and let $\mathbf x=(h+1)\mathbf e_i$.
Then
\[
\sum_{j=1}^k x_j=h+1,
\qquad
\mathbf x\cdot\vec{A}=(h+1)a=n.
\]
If $b>0$, then $a<q$. Choose $A\in\binom{[q]}k$ containing
$a=a_i$ and $a+1=a_j$, and let
$\mathbf x=(h+1-b)\mathbf e_i+b\mathbf e_j$. Then
\begin{align*}
\sum_{t=1}^k x_t
& \ = \ (h+1-b)+b \ = \ h+1,\\
\mathbf x\cdot\vec{A}
& \ = \ (h+1-b)a+b(a+1)\\
& \ = \ (h+1)a+b \ = \ n.
\end{align*}
These choices of $A$ are possible because $q\geq k$. Therefore
\begin{align*}
\mathcal S_{h+1,k}(q)
& \ = \ [(h+1)..(h+1)q],\\
\abs{\mathcal S_{h+1,k}(q)}
& \ = \ (h+1)q-(h+1)+1 \ = \ (h+1)q-h.
\end{align*}

For the symmetry, define
\[
A'=\{q+1-a:a\in A\},
\qquad
\mathbf x'=(x_k,\ldots,x_1),
\qquad
\vec{A'}=(q+1-a_k,\ldots,q+1-a_1).
\]
Then
\begin{align*}
\mathbf x'\cdot\vec{A'}
& \ = \ \sum_{i=1}^k x_{k+1-i}(q+1-a_{k+1-i})\\
& \ = \ (q+1)\sum_{i=1}^k x_i-\sum_{i=1}^k x_i a_i\\
& \ = \ (h+1)(q+1)-\mathbf x\cdot\vec{A}.
\end{align*}
Applying the same transformation twice gives $A''=A$ and
$\mathbf x''=\mathbf x$. It is therefore a bijection between the pairs
counted by $\nu_q(s)$ and those counted by
$\nu_q((h+1)(q+1)-s)$, giving the stated equality.
\end{proof}

\subsection{Questions}
We conclude with several questions suggested by the results above.

\begin{question}\label{q:threshold}
For $h\geq2$ and $k\geq3$, let $q_0(h,k)$ be the smallest integer such that Corollary~\ref{cor:gap-pattern} holds for every $q\geq q_0(h,k)$. By Remark~\ref{rem:explicit},
\[
q_0(h,k)\leq 2^k k!h(2h+1)^{2k-2}+1.
\]
Can this bound be improved? How fast can $h$ grow with $q$ while the conclusion still holds?

Direct computation gives the following possible values of $q_0(h,k)$:
\[
\begin{array}{c|cccc}
(h,k)&(3,3)&(4,3)&(2,4)&(3,4)\\
\hline
q_0(h,k)&5&6&8&16
\end{array}
\]
These values have been verified for $q\leq120,120,60,50$, respectively, with the conclusion failing at $q=4,5,7,15$.
\end{question}

\begin{question}\label{q:rate}
By Corollary~\ref{cor:proportion},
\[
q\left(1-\frac{\abs{\Bcal_{h,k}(q)}}{\binom qk}\right)
\]
is bounded above and below by positive constants. For $k=3$, it converges to $3\sum_{w=2}^h\varphi(w)/w$ by Remark~\ref{rem:k3-limit}. Does it converge for every $k\geq4$? If so, what is its limit $c_{h,k}$? Is it eventually monotone in $q$?

For $(h,k)=(2,4)$, the only relations of weight at most $2$ that can hold when $a_1<a_2<a_3<a_4$ are
\[
2a_2=a_1+a_3,\qquad 2a_3=a_2+a_4,\qquad 2a_2=a_1+a_4,\qquad 2a_3=a_1+a_4,\qquad a_1+a_4=a_2+a_3.
\]
The numbers of sets in $\binom{[q]}4$ satisfying them are asymptotic to $q^3/12$, $q^3/12$, $q^3/24$, $q^3/24$ and $q^3/12$, respectively, and only $O(q^2)$ sets satisfy two of them. Hence $c_{2,4}=24\cdot\frac13=8$.

Can $c_{h,k}$ be expressed in general as a sum of normalized $(k-1)$-dimensional volumes of the sections
\[
\left\{
\mathbf v\in\R^k:
0<v_1<\cdots<v_k<1,\
\mathbf d\cdot\mathbf v=0
\right\},
\qquad \mathbf d\in\Rcal_{h,k}?
\]
\end{question}

\begin{question}\label{q:near-collisions}
For an integer $t$, let $N_{h,k}(t;q)$ count the sets $A\in\binom{[q]}k$ for which there exist distinct $\mathbf x,\mathbf y\in\X_{h,k}$ satisfying
\[
(\mathbf x-\mathbf y)\cdot\vec A=t.
\]
In particular, $N_{h,k}(0;q)$ counts the sets that are not $B_h$-sets. Each of the at most $M_{h,k}^2$ pairs $(\mathbf x,\mathbf y)$ gives a single equation $(\mathbf x-\mathbf y)\cdot\mathbf v=t$ with $\mathbf x-\mathbf y\neq\mathbf 0$, so by the argument of Lemma~\ref{lem:rank-lemma},
\[
N_{h,k}(t;q) \ \leq\ M_{h,k}^2q^{k-1}.
\]
How does $N_{h,k}(t;q)$ depend on $t$? More generally, what can be said about the distribution of differences in $hA-hA$?

Related questions on sums and differences of random sets have been studied in~\cite{HegartyMiller2009,ILMZ2012,MartinOBryant2007}.
\end{question}

\begin{question}\label{q:distribution}
Can the function $\nu_q$ from Proposition~\ref{prop:S-interval} be determined explicitly? In particular, does
\[
\lim_{q\to\infty}
q^{-(k-1)}\nu_q(\lfloor\sigma q\rfloor)
\]
exist for every $\sigma\in(0,h+1)$? If so, is the resulting function piecewise polynomial?
\end{question}

\bibliographystyle{amsalpha}
\bibliography{references}

\end{document}